\documentclass{article}
\usepackage[centertags]{amsmath}
\usepackage{amsfonts}
\usepackage{amssymb}
\usepackage{amsthm}
\usepackage{newlfont}
\usepackage{graphicx}
\newcommand \qq {\qquad}
\newcommand \q  {\quad}

\newfont{\bb}{msbm10 at 10pt}
\def\r{\hbox{\bb R}}

\theoremstyle{plain}
\newtheorem{thm}{Theorem}[section]

\theoremstyle{definition}

\numberwithin{equation}{section}

\begin{document}

\centerline{\textbf{ \begin{large} Constant Scalar Curvature of
Three Dimensional Surfaces
\end{large}}}
\centerline{\textbf{ \begin{large}  Obtained by the Equiform Motion
of a helix
\end{large}}}

\bigskip
\centerline{\begin{large} Ahmad T. Ali $^{a}$\footnote{Corresponding author.\\
E-mail address: atali71@yahoo.com (A. T. Ali).}, Fathi M.
Hamdoon$^b$ and Rafael L\'{o}pez$^c$
\end{large} }

\centerline{\footnotesize $^a$ Mathematics Department, Faculty of
Science, Al-Azhar University, Cairo, Egypt
}\centerline{\footnotesize $^b$ Mathematics Department, Faculty of
Science, Al-Azhar University, Assiut, Egypt }
\centerline{\footnotesize $^c$ Departamento de Geometr\'{\i}a y
Topolog\'{\i}a, Universidad de Granada, 18071 Granada, Spain } $ $
\hrule

\begin{flushleft}\textbf{Abstract}.\end{flushleft}

{\it In this paper we consider the equiform motion of a helix in
Euclidean space $\mathbf{E}^7$. We study and analyze the
corresponding kinematic three dimensional surface under the
hypothesis that its scalar curvature $\mathbf{K}$ is constant. Under
this assumption, we prove that if the scalar curvature $\mathbf{K}$
is constant, then $\mathbf{K}$ must equal zero.}

\begin{flushleft}
\end{flushleft}

\begin{flushleft}
\end{flushleft}
\smallskip
\smallskip
\hrule

\section{Introduction}

Helix is one of the most fascinating curves in science and nature. Scientist have long held a fascinating, sometimes bordering on mystical obsession, for helical structures in nature. Helices arise in nano-springs, carbon nano-tubes, $\alpha$-helices, DNA double and collagen triple helix, lipid bilayers, bacterial flagella in salmonella and escherichia coli, aerial hyphae in actinomycetes, bacterial shape in spirochetes, horns, tendrils, vines, screws, springs, helical staircases and sea shells \cite{choua, lucas, watson}. Also we can see the helix curve or helical structures in fractal geometry, for instance hyperhelices \cite{toledo}. In the field of computer aided design and computer graphics, helices can be used for the tool path description, the simulation of kinematic motion or the design of highways, etc. \cite{yang}. From the view of differential geometry, a helix is a geometric curve with non-vanishing constant curvature $\kappa$ and non-vanishing constant torsion $\tau$ \cite{barros}.

An equiform transformation in the $n$-dimensional Euclidean space
$\mathbf{E}^n$ is an affine transformation whose linear part is
composed by an orthogonal transformation and a homothetical
transformation. Such an equiform transformation maps points
$\mathbf{x}\in\mathbf{E}^n$ according to the rule
\begin{equation}\label{equi1}
\mathbf{x}\longmapsto s {\cal A}\mathbf{x}+
\mathbf{d},\hspace*{.5cm} {\cal A}\in SO(n), s\in\mathbf{R}^+,
\mathbf{d}\in\mathbf{E}^n.
\end{equation}
The number $s$ is called the scaling factor. An equiform motion is
defined if the parameters of (\ref{equi1}), including $s$, are given
as functions of a time parameter $t$. Then a smooth one-parameter
equiform motion moves a point $\mathbf{x}$ via
$\mathbf{x}(t)=s(t){\cal A}(t) \mathbf{x}(t)+\mathbf{d}(t)$. The
kinematic corresponding to this transformation group is called
equiform kinematic. See  \cite{br, fhk}.

Under the assumption of the constancy of the scalar curvature,
kinematic surfaces obtained by the motion of a circle have been obtained in \cite{ah}. In
a similar context, one can consider hypersurfaces in  space forms
generated by  one-parameter family of spheres and having constant
 curvature: \cite{ca, ja, lo, pa}.

In this paper we consider the equiform motions of a helix
$\mathbf{k}_0$ in $\mathbf{E}^n$. The point paths of the helix
generate a $3$-dimensional surface, containing  the positions of the
starting helix $\mathbf{k}_0$. We have studied the first order
properties of these surfaces for arbitrary dimensions $n\geq 3$ and
for the points of these helixes. We restrict our considerations to
dimension $n=7$ because at any moment, the infinitesimal
transformations of the motion map the points of the helix
$\mathbf{k}_0$ to the velocity vectors, whose end points will form
an affine image of $\mathbf{k}_0$ (in general a helix) that span a
subspace $\mathbf{W}$ of $\mathbf{E}^n$ with $n\leq7$.

Let $\mathbf{x}(\phi)$ be a parameterization of $\mathbf{k}_0$ and
let $\mathbf{X}(t,\phi)$ be the resultant 3-surface by the
equiform motion. We consider a certain position of the moving
space given by $t=0$, and we would like to obtain information
about the motion at  least during a certain period around $t=0$ if
we know its characteristics for one instant. Then we  restrict our
study to the properties of the motion for the limit case
$t\rightarrow 0$. A first choice is then approximate
$\mathbf{X}(t,\phi)$ by the first derivative of the trajectories.
Solliman, et al. studied $3$-dimensional surfaces in
$\mathbf{E}^7$ generated by equiform motions of a sphere proving
that, in general, they are contained in a canal hypersurface
\cite{skhs}.

The purpose of this paper is to describe the kinematic surfaces
obtained by the motion of a helix and whose scalar curvature
$\mathbf{K}$ is constant. As a consequence of our results, we
prove:

\begin{quote} \emph{ Consider a kinematic three-dimensional surface obtained by
the equiform motion of a helix . If the scalar curvature
$\mathbf{K}$ is constant, then $\mathbf{K}$ must vanish on the
surface.}
\end{quote}
Moreover, and in the case that $\mathbf{K}=0$ we show the
description of  the motion of such 3-surface giving the equations
that determine the kinematic geometry. We will show an example of a
such surface.

\section{The representation of a kinematic surface}

In two copies $\sum^{0}$, $\sum$ of Euclidean 7-space
$\mathbf{E}^7$,  we consider a circular helix  $\mathbf{k}_0$ at a
cylinder with unit radius centered at the origin of the 3-space
$\varepsilon_0=[x_1x_2x_3]$ and represented by
\begin{equation}\label{eq21}\mathbf{x}(\phi)=\Big(\cos(\phi),
\sin(\phi),\lambda\,\phi,0,0,0,0\Big)^{\text{T}},\q\,\,\,
\phi\in[0,2 \pi].
\end{equation}
Under a one-parameter equiform motion of moving space $\sum^{0}$
with respect to a fixed space $\sum$ the general representation of
the motion of this surface  in $\mathbf{E}^7$ is given by
$$\mathbf{X}(t,\phi)=s(t)\mathbf{A}(t)\mathbf{x}(\phi)+\mathbf{d}(t),\qq t\in I\subset
\mathbf{R}.$$
 Here
$\mathbf{d}(t)=\Big(b_i(t)\Big)^{\text{T}}:i=1,2,...,7$ describes
the position of the origin of $\sum^{0}$ at time $t$,
$\mathbf{A}(t)=\Big(a_{ij}(t)\Big)^{\text{T}}:i,j=1,2,...,7$ is an
orthogonal matrix and $s(t)$ provides the scaling factor of the
moving system. With $s=\text{const.}\neq 0$ (sufficient to set
$s=1$), we have an ordinary Euclidean rigid body motion. For
varying $t$ and fixed $\mathbf{x}(\phi)$, equation (\ref{eq21})
gives a parametric representation of the surface (or trajectory)
of $\mathbf{x}(\phi)$. Moreover, we assume that all involved
functions are at least of class $\mathbf{C}^1$. Using Taylor's
expansion up to the first order,  the representation of the motion
is given by
$$\mathbf{X}(t,\phi)=\Big[s(0)\mathbf{A}(0)+t\Big(\dot{s}(0)\mathbf{A}(0)+s(0)\dot{\mathbf{A}}(0)
\Big)\Big]\mathbf{x}(\phi)+\mathbf{d}(0)+t\,\dot{\mathbf{d}}(0),
$$ where $(.)$ denotes differentiation with respect to the time $t$.
Assuming that the moving frames $\sum^{0}$ and $\sum$ coincide at
the zero position ($t=0$),  we have
$$\mathbf{A}(0)=\mathbf{I},\q\, s(0)=1\q\, \text{and}\q\, \mathbf{d}(0)=0.
$$
Thus we have
$$\mathbf{X}(t,\phi)=\Big[\mathbf{I}+t\Big(s^{\prime}\mathbf{I}+\Omega
\Big)\Big]\mathbf{x}(\phi)+t\mathbf{d}^{\prime},$$
 where $\Omega=\dot{\mathbf{A}}(0)=(\omega_i),\,i=1,2,...,21$ is a skew symmetric matrix,
$s^{\prime}=\dot{s}(0)$, $\mathbf{d}^{\prime}=\dot{\mathbf{d}}(0)$
and all values of $s, b_i$ and their derivatives are computed at
$t=0$. With respect to these frames, the representation of the
motion up to the first order is
\begin{eqnarray}
\begin{pmatrix}
  \mathbf{X}_1 \\
   \mathbf{X}_2 \\
  \mathbf{X}_3 \\
   \mathbf{X}_4 \\
   \mathbf{X}_5 \\
   \mathbf{X}_6 \\
   \mathbf{X}_7
\end{pmatrix}=\begin{pmatrix}
 1+s^{\prime}t & \omega_1\,t & \omega_2\,t & \omega_3\,t & \omega_4\,t & \omega_5\,t & \omega_6\,t \\
  -\omega_1\,t & 1+s^{\prime}t & \omega_7\,t & \omega_8\,t & \omega_9\,t & \omega_{10}\,t & \omega_{11}\,t \\
  -\omega_2\,t & -\omega_7\,t & 1+s^{\prime}t & \omega_{12}\,t & \omega_{13}\,t & \omega_{14}\,t & \omega_{15}\,t \\
  -\omega_3\,t & -\omega_8\,t & -\omega_{12}\,t & 1+s^{\prime}t & \omega_{16}\,t & \omega_{17}\,t & \omega_{18}\,t \\
  -\omega_4\,t & -\omega_9\,t & -\omega_{13}\,t & -\omega_{16}\,t & 1+s^{\prime}t & \omega_{19}\,t & \omega_{20}\,t \\
  -\omega_5\,t & \omega_{10}\,t & -\omega_{14}\,t & -\omega_{17}\,t & -\omega_{19}\,t & 1+s^{\prime}t & \omega_{21}\,t \\
  -\omega_6\,t & -\omega_{11}\,t & -\omega_{15}\,t & -\omega_{18}\,t & -\omega_{20}\,t & -\omega_{21}\,t & 1+s^{\prime}t
\end{pmatrix}\,\begin{pmatrix}
 \cos(\phi) \\
 \sin(\phi) \\
 \lambda\,\phi \\
  0 \\
  0 \\
  0 \\
  0
\end{pmatrix}+t\begin{pmatrix}
  b^{\prime}_{1} \\
  b^{\prime}_{2} \\
  b^{\prime}_{3} \\
  b^{\prime}_{4} \\
  b^{\prime}_{5} \\
  b^{\prime}_{6} \\
  b^{\prime}_{7}
\end{pmatrix},\nonumber
\end{eqnarray}
or in the equivalent form
\begin{eqnarray}\begin{pmatrix}
  \mathbf{X}_{1} \\
  \mathbf{X}_{2} \\
  \mathbf{X}_{3} \\
  \mathbf{X}_{4} \\
  \mathbf{X}_{5} \\
  \mathbf{X}_{6} \\
  \mathbf{X}_{7}
\end{pmatrix}=t\begin{pmatrix}
  b^{\prime}_{1} \\
  b^{\prime}_{2} \\
  b^{\prime}_{3} \\
  b^{\prime}_{4} \\
  b^{\prime}_{5} \\
  b^{\prime}_{6} \\
  b^{\prime}_{7}
\end{pmatrix}+\cos(\phi)\begin{pmatrix}
  1+s^{\prime}\,t \\
  -\omega_{1}\,t \\
  -\omega_{2}\,t \\
  -\omega_{3}\,t \\
  -\omega_{4}\,t \\
  -\omega_{5}\,t \\
  -\omega_{6}\,t
\end{pmatrix}+\sin(\phi)\begin{pmatrix}
  \omega_{1}\,t \\
  1+s^{\prime}\,t \\
  -\omega_{7}\,t \\
  -\omega_{8}\,t \\
  -\omega_{9}\,t \\
  -\omega_{10}\,t \\
  -\omega_{11}\,t
\end{pmatrix}+\lambda\,\phi\begin{pmatrix}
  \omega_{2}\,t \\
   \omega_{7}\,t \\
  1+s^{\prime}\,t \\
  -\omega_{12}\,t \\
  -\omega_{13}\,t \\
  -\omega_{14}\,t \\
  -\omega_{15}\,t
\end{pmatrix}.\label{eq22}
\end{eqnarray}

\section{ Scalar curvature of the kinematic surface }

In this section we shall compute the scalar curvature of 3-surfaces
in $\mathbf{E}^7$ generated by equiform motions of a helix. The
tangents to the parametric curves $t=\text{const}.$ and
$\phi=\text{const}.$ at the zero position  are
$$\mathbf{X}_t=\Big[s^{\prime}\mathbf{I}+\Omega
\Big]\mathbf{x}+\mathbf{d}^{\prime},\qq\,\,
\mathbf{X}_\phi=\Big[\mathbf{I}+\Big(s^{\prime}\mathbf{I}+\Omega\Big)t
\Big]\mathbf{x}_{\phi},
$$
The first fundamental quantities of $\mathbf{X}(t,\phi)$ are
\begin{eqnarray*}
g_{11}=\mathbf{X}_t\,\mathbf{X}_t^{\text{T}},\q\,
g_{12}=\mathbf{X}_\phi\,\mathbf{X}_t^{\text{T}},\q\,
g_{22}=\mathbf{X}_\phi\,\mathbf{X}_\phi^{\text{T}}.
\end{eqnarray*}
Now, we obtain
\begin{eqnarray*}
g_{11}&=&\alpha_{0}+\alpha_{1}\,\phi+\alpha_{2}\,\phi^2+\Big(\alpha_{3}+\alpha_{4}\,\phi\Big)\cos(\phi)
+\Big(\alpha_{5}+\alpha_{6}\,\phi\Big)\sin(\phi)+\alpha_{7}\cos(2\phi)+\alpha_{8}\sin(2\phi),\\
g_{12}&=&\lambda\,b_3^{\prime}-\omega_1+\lambda^2\,s^{\prime}\,\phi+\Big(b_2^{\prime}-\lambda\,\omega_2
+\lambda\,\omega_7\,\phi\Big)\cos(\phi) -\Big(
b_1^{\prime}+\lambda\,\omega_7+\lambda\,\omega_2\,\phi
\Big)\sin(\phi)+\dfrac{1}{2}\,t\,\Big[ \alpha_{1}\\
&+&2\alpha_{2}\,\phi
+\Big(\alpha_{4}+\alpha_{5}+\alpha_{6}\phi\Big)\cos(\phi)
+\Big(\alpha_{6}-\alpha_{3}-\alpha_{4}\phi\Big)\sin(\phi)
+2\alpha_{8}\cos(2\phi)-2\alpha_{7}\sin(2\phi)
\Big],\\
g_{22}&=&\Big(1+\lambda^2\Big)\Big(1+2\,s^{\prime}\,t\Big)
+t^2\Big[\alpha_{9}+\alpha_{6}\cos(\phi)-\alpha_{7}\cos(2\phi)-\alpha_{4}\sin(\phi)
-\alpha_{8}\sin(2\phi)\Big].
\end{eqnarray*}
where
\begin{eqnarray*}
\begin{cases}
\alpha_{0}=s^{\prime2}+\omega_1^2+\dfrac{1}{2}\Big[\sum_{i=1}^7
b_i^2+\sum_{i=2}^{11}\omega_{i}^2\Big], & \\
\alpha_{1}=2\,\lambda\,\Big[b_1^{\prime}\,\omega_2+b_2^{\prime}\,\omega_2+b_3^{\prime}\,s^{\prime}-\sum_{i=4}^{7}
b_{i}^{\prime}\,\omega_{i+8}\Big], & \\
\alpha_{2}=\lambda^2\,\Big[s^{\prime2}+\omega_2^2+\omega_7^2+\sum_{i=12}^{15}
\omega_{i}^2\Big],& \\ \alpha_{3}=2\,\Big[b_1^{\prime}\,s^{\prime}-\sum_{i=1}^{6}b_{i+1}^{\prime}\omega_{i}\Big], & \\
\alpha_{4}=
2\,\lambda\,\Big[-\omega_1\,\omega_7+\sum_{i=3}^{6}\omega_{i}\,\omega_{i+9}\Big],
& \\
\alpha_{5}=2\,\Big[b^{\prime}_1\,\omega_1+b^{\prime}_2\,s^{\prime}-\sum_{i=3}^{7}
b^{\prime}_{i}\,\omega_{i+4}\Big], & \\
\alpha_{6}=2\,\lambda\,\Big[\omega_1\,\omega_2+\sum_{i=8}^{11}
\omega_{i}\,\omega_{i+4}\Big], & \\
\alpha_{7}=\dfrac{1}{2}\,\sum_{i=2}^{6}\Big(
\omega_{i}^2-\omega_{i+5}^2\Big), & \\
\alpha_{8}=\sum_{i=2}^{6} \omega_{i}\,\omega_{i+5}, & \\
\alpha_{9}=(1+\lambda^2)\,s^{\prime2}+\omega_1^2+\dfrac{1}{2}\,
\sum_{i=2}^{11}\omega_i^2+\lambda^2\Big(
\omega_2^2+\omega_2^7+\sum_{i=12}^{15}\omega_i^2\Big).
\end{cases}
\end{eqnarray*}
In order to calculate the scalar curvature, we need to compute the
Christoffel symbols of the second kind, which are defined as
\begin{equation}\label{eq311}
\Gamma_{ij}^{l}=\dfrac{1}{2}g^{lm}\Big[\frac{\partial
g_{im}}{\partial x_{j}}+\frac{\partial g_{jm}}{\partial
x_{i}}-\frac{\partial g_{ij}}{\partial x_{m}}\Big],
\end{equation}
where $i,j,l$ are indices that take the values $1,2$, $x_1=t,
x_2=\phi$, and $\Big(g^{lm}\Big)$ is the inverse matrix of
$\Big(g_{ij}\Big)$. Then the scalar curvature of the surface
$\mathbf{X}(t,\phi)$ is
$$\mathbf{K}(t,\phi)=g^{ij}\Big[\frac{\partial
\Gamma_{ij}^l}{\partial x_{l}}-\frac{\partial
\Gamma_{il}^l}{\partial
x_{j}}+\Gamma_{ij}^l\,\Gamma_{lm}^m-\Gamma_{il}^m\,\Gamma_{jm}^l\Big].$$
At the zero position ($t=0$), the scalar curvature of
$\mathbf{X}(t,\phi)$ is given by
\begin{eqnarray}
\mathbf{K}=\mathbf{K}(0,\phi)&=&\dfrac{P\Big(\phi^{n_1}\cos(m_1\,\phi),\phi^{n_1}\sin(m_1\,\phi)\Big)}{
Q\Big(\phi^{n_2}\cos(m_2\,\phi),\phi^{n_2}\sin(m_2\,\phi)\Big)}\nonumber\\
&=&
\dfrac{\sum_{j=0}^{6}\sum_{i=0}^{4}\Big(A_{i,j}\phi^{i}\cos(j\,\phi)
+B_{i,j}\phi^i\sin(j\,\phi)\Big)}{\sum_{j=0}^{6}\sum_{i=0}^{6}\Big(F_{i,j}\phi^{i}\cos(j\,\phi)
+H_{i,j}\phi^i\sin(j\,\phi)\Big)}.\label{eq31}
\end{eqnarray}
This quotient writes then as
\begin{equation}\label{eq32}
P\Big(\phi^{n_1}\cos(m_1\,\phi),\phi^{n_1}\sin(m_1\,\phi)\Big)
-\mathbf{K}\,Q\Big(\phi^{n_2}\cos(m_2\,\phi),\phi^{n_2}\sin(m_2\,\phi)\Big)=0.
\end{equation}
The assumption on the constancy of the scalar curvature
$\mathbf{K}$ implies that  equation (\ref{eq32}) is a linear
combination of the functions $\phi^i$, $\cos(n\,\phi)$ and
$\sin(m\,\phi)$. Because these functions are independent linearly,
the corresponding coefficients must  vanish. Throughout this work,
we have employed the Mathematica programm in order to compute the
explicit expressions of these coefficients.

In the next two subsections, we distinguish the cases $\mathbf{K}=0$
and $\mathbf{K}\not=0$ respectively.

\subsection{ Kinematic surfaces with zero scalar curvature}

We assume that $\mathbf{K}=0$. From the expression (\ref{eq31}), we
have
\begin{eqnarray*}
P\Big(\phi^{n}\cos(m\,\phi),\phi^{n}\sin(m\,\phi)\Big)=\sum_{j=0}^{6}\sum_{i=0}^{4}\Big(A_{i,j}\phi^{i}\cos(j\,\phi)
+B_{i,j}\phi^i\sin(j\,\phi)\Big)=0.
\end{eqnarray*}
In this case, a straightforward computation shows that the
coefficients of $\phi^4\cos(6\phi)$ and $\phi^4\sin(6\phi)$ are
\begin{eqnarray*}
A_{4,6}&=&\frac{1}{4}\,\lambda^4\,\Big[4\,\alpha_{8}\,\omega_2\,\omega_7\,\Big(\omega_2^2-\omega_7^2\Big)
-\alpha_{7}\Big(\omega_2^4-6\,\omega_2^2\,\omega_4^2+\omega_7^2\Big)\Big], \nonumber\\
B_{4,6}&=&\frac{1}{4}\,\lambda^4\,\Big[4\,\alpha_{7}\,\omega_2\,\omega_7\,\Big(\omega_7^2-\omega_2^2\Big)
-\alpha_{8}\Big(\omega_2^4-6\,\omega_2^2\,\omega_4^2+\omega_7^2\Big)\Big].
\nonumber
\end{eqnarray*}
We consider the two equations $A_{4,6}=0$ and $B_{4,6}=0$. From
here, we discuss two possibilities: $\alpha_{7}=\alpha_{8}=0$ and
$\omega_2=\omega_7=0$.

\begin{enumerate}
\item {\bf Case} $\alpha_{7}=\alpha_{8}=0$.

A computation of coefficients yields
\begin{eqnarray*}
A_{4,5}&=&\frac{1}{4}\,\lambda^4\,\Big[4\,\alpha_{4}\,\omega_2\,\omega_7\,\Big(\omega_2^2-\omega_7^2\Big)
+\alpha_{6}\Big(\omega_2^4-6\,\omega_2^2\,\omega_4^2+\omega_7^2\Big)\Big]=0, \nonumber\\
B_{4,5}&=&\frac{1}{4}\,\lambda^4\,\Big[4\,\alpha_{6}\,\omega_2\,\omega_7\,\Big(\omega_2^2-\omega_7^2\Big)
-\alpha_{4}\Big(\omega_2^4-6\,\omega_2^2\,\omega_4^2+\omega_7^2\Big)\Big]=0.
\nonumber
\end{eqnarray*}
We consider two cases: $\alpha_{4}=\alpha_{6}=0$ and
$\omega_2=\omega_7=0$.

\begin{enumerate}
\item  We assume $\alpha_{4}=\alpha_{6}=0$. The computation of
coefficients leads to
\begin{eqnarray*}
A_{4,4}&=&\frac{1}{2}\,\lambda^4\,\alpha_{9}\,\Big(\omega_2^4-6\,\omega_2^2\,\omega_4^2+\omega_7^2\Big)=0, \nonumber\\
B_{4,4}&=&2\,\lambda^4\,\omega_2\,\omega_7\,\alpha_{9}\Big(\omega_2^2-\omega_7^2\Big)=0.
\nonumber
\end{eqnarray*}
Because $\alpha_{9}\neq0$, we conclude $\omega_2=\omega_7=0$. New
computations give
\begin{eqnarray*}
A_{4,0}&=&4\Big[(1+\lambda^2)\alpha_{2}-\lambda^4
s^{\prime2}\Big]\Big[
(\lambda^2+1)\alpha_{2}(\alpha_{9}-\alpha_{2})+s^{\prime2}\Big((\lambda^4-1)\alpha_{2}-\lambda^4
\alpha_{9}\Big)\Big]. \nonumber
\end{eqnarray*}
The first term in $A_{4,0}$ is not zero because
\begin{equation}\label{a02}
(1+\lambda^2)\alpha_{2}-\lambda^4\,s^{\prime2}=\lambda^2\Big[s^{\prime2}+(1+\lambda^2)\sum_{i=12}^{15}\omega_i^2\Big]>0.
\end{equation}
The second term take the following form:
$$
\lambda^2 s^{\prime2}
\Big[\omega_1^2+\frac{1}{2}\sum_{i=3}^{6}\Big(\omega_i^2+\omega_{i+5}^2\Big)
+\lambda^2
\sum_{i=12}^{15}\omega_i^2\Big]+(\lambda^4+\lambda^2)\Big[\omega_1^2+\frac{1}{2}\sum_{i=3}^{6}\Big(\omega_i^2
+\omega_{i+5}^2\Big)\Big]\Big[\sum_{i=12}^{15}\omega_i^2\Big].
$$
When $A_{4,0}=0$ must leads the following conditions:
$$\omega_{k}=0,\,k=1,3,4,5,6,8,9,...,15. $$
From here, all coefficients $A_{i,j}$ and $B_{i,j}$ are equal zero.

\item  We now assume that $\omega_2=\omega_7=0$ and that
$\alpha_{4},\,\alpha_{6}\neq0$. The computation of coefficients
yields
$$
A_{4,0}=4\,\lambda^2\,s^{\prime2}\alpha_{6}\Big((1+\lambda^2)\alpha_{2}-\lambda^4\,s^{\prime2}\Big)=0.
$$
Then $\alpha_{6}=0$ and hence, we obtain
$$
B_{4,0}=4\,\lambda^2\,s^{\prime2}\alpha_{4}\Big(\lambda^4\,s^{\prime2}-(1+\lambda^2)\alpha_{2}\Big)=0.
$$
which leads $\alpha_{4}=0$ contradiction.
\end{enumerate}

\item {\bf Case}  $\omega_2=\omega_7=0$ and
$\alpha_{7},\,\alpha_{8}\neq0$.

The computation of coefficients yields
$$
A_{4,1}=4\,\lambda^2\,s^{\prime2}\alpha_{6}\Big((1+\lambda^2)\alpha_{2}-\lambda^4\,s^{\prime2}\Big)$$
$$B_{4,1}=4\,\lambda^2\,s^{\prime2}\alpha_{4}\Big(\lambda^4\,s^{\prime2}-(1+\lambda^2)\alpha_{2}\Big).
$$
Then $\alpha_{6}=\alpha_{4}=0$, the new computations give
$$
A_{4,2}=-4\,\alpha_{7}\Big((1+\lambda^2)\alpha_{2}-\lambda^4\,s^{\prime2}\Big)^2$$
$$B_{4,2}=-4\,\alpha_{8}\Big((1+\lambda^2)\alpha_{2}-\lambda^4\,s^{\prime2}\Big)^2.
$$
By solving the equations $A_{4,2}=0$ and $B_{4,2}=0$, we get
$\alpha_{7}=\alpha_{8}=0$: contradiction.
\end{enumerate}
As a consequence of the above computations, the only case that
occurs is when $\alpha_{4}=\alpha_{6}=\alpha_{7}=\alpha_{8}=0$.
This can summarize as follows:

\begin{thm}\label{th-31} A kinematic 3-surface in $\mathbf{E}^7$ generated by the equifom motion of a helix,
with zero scalar curvature and whose parametrization writes as
(\ref{eq22}), satisfies
$$\omega_{k}=0,\hspace*{1cm}1\leq k\leq 15.$$
\end{thm}

In order to end with this subsection, we show an example of a such
$3$-surface. Let us consider
\begin{eqnarray}
\bold{A}(t)=\begin{pmatrix}
 \cos t & 0 & 0 & 0 & 0 & \sin t\sin\mu t & 0 \\
  0 & \cos\mu t & 0 & 0 & \sin t\sin\mu t & 0 & 0 \\
  0 & 0 & \cos\mu t & 0 & 0 & 0 & -\sin t\sin\mu t \\
  0 & 0 & 0 & \cos\mu t & \sin\mu t & 0 & 0 \\
  0 & 0 & 0 & -\sin\mu t & \cos\mu t & 0 & \sin\mu t\cos\mu t \\
  -\sin t\sin\mu\,t & 0 & 0 & 0 & 0 & \cos\mu t & \sin\mu t \\
  \sin t\sin\mu t & 0 & 0 & 0 & -\sin\mu t & -\cos t\sin\mu t &
  \cos\mu t
\end{pmatrix}.\nonumber
\end{eqnarray}
such that $\mu\in \r-\{0\}$. We assume $s(t)=e^{s^{\prime}\,t}$ and
$ \bold{d}(t)=t\,\Big( b_1^{\prime}, b_2^{\prime}, b_3^{\prime},
b_4^{\prime}, b_5^{\prime}, b_6^{\prime}, b_7^{\prime}\Big)^{T}$. By
differentiating $\bold{A}(t)$ and putting $t=0$, we find
$$
\omega_{16}=\omega_{20}=\omega_{21}=\mu,\,\,\,\text{and}\,\,\,\omega_{k}=0,\,\,\,k=1,2,...,15,17,18,19.
$$
From Theorem \ref{th-31}, and by taking arbitrary numbers
$\lambda$, $\mu$ and $b_i^{\prime}$, $1\leq i\leq 7$, the
corresponding kinematic $3$-surface given by Equation (\ref{eq22})
has zero constant scalar curvature. We show a picture of a such
surface
 $\bold{X}(t,\phi)$ in axonometric viewpoint $\bold{Y}(t,\phi)$. For
this, the unit vectors $E_4=(0,0,0,1,0,0,0)$, $E_5=(0,0,0,0,1,0,0)$,
$E_6=(0,0,0,0,0,1,0)$ and $E_7=(0,0,0,0,0,0,1)$ are mapped onto the
vectors $(1,1,0)$, $(1,0,1)$, $(0,1,1)$ and $(1,1,1)$ respectively
(see  \cite{gss}). With the above choice,
\begin{eqnarray}
  \mathbf{X}(t,\phi)=t\begin{pmatrix}
  b^{\prime}_{1} \\
  b^{\prime}_{2} \\
  b^{\prime}_{3} \\
  b^{\prime}_{4} \\
  b^{\prime}_{5} \\
  b^{\prime}_{6} \\
  b^{\prime}_{7}
\end{pmatrix}+\cos(\phi)\begin{pmatrix}
  1+s^{\prime}\,t \\
  0 \\
  0 \\
  0 \\
  0 \\
  0 \\
  0
\end{pmatrix}+\sin(\phi)\begin{pmatrix}
  0 \\
  1+s^{\prime}\,t \\
  0 \\
  0 \\
  0 \\
  0 \\
  0
\end{pmatrix}+\lambda\,\phi\begin{pmatrix}
  0 \\
   0 \\
  1+s^{\prime}\,t \\
  0 \\
  0 \\
  0 \\
  0
\end{pmatrix}.\nonumber
\end{eqnarray}
and
\begin{eqnarray}
  \mathbf{Y}(t,\phi)=t\begin{pmatrix}
  B^{\prime}_{1} \\
  B^{\prime}_{2} \\
  B^{\prime}_{3}
\end{pmatrix}+\cos(\phi)\begin{pmatrix}
  1+s^{\prime}\,t \\
  0 \\
  0
\end{pmatrix}+\sin(\phi)\begin{pmatrix}
  0 \\
  1+s^{\prime}\,t \\
  0
\end{pmatrix}+\lambda\,\phi\begin{pmatrix}
  0 \\
   0 \\
  1+s^{\prime}\,t
\end{pmatrix},\nonumber
\end{eqnarray}
where
$B_1^{\prime}=b_1^{\prime}+b_4^{\prime}+b_5^{\prime}+b_7^{\prime}$,
$B_2^{\prime}=b_2^{\prime}+b_4^{\prime}+b_6^{\prime}+b_7^{\prime}$
and
$B_3^{\prime}=b_3^{\prime}+b_5^{\prime}+b_6^{\prime}+b_7^{\prime}$.
Let us fix the following values:
$$
\lambda=\frac{1}{2\pi},\hspace*{1cm}B_1=B_2=B_3=s'=1.
$$

\begin{figure}[h]
\includegraphics[width=6cm]{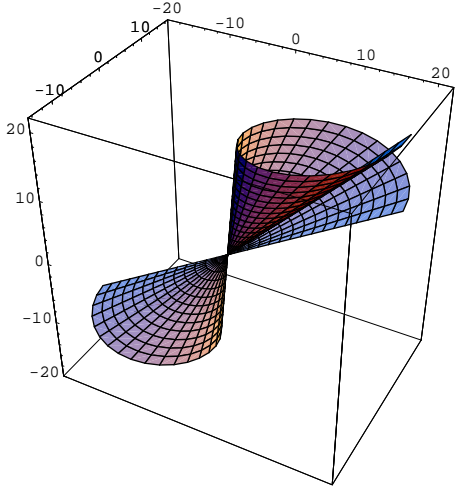}
\caption{A kinematic 3-surface given by the motion of a helix and
with $\mathbf{K}=0$. }\label{fig1}
\end{figure}

The corresponding picture is given in Figure \ref{fig1}.

\subsection{ Kinematic surfaces with non-zero constant scalar curvature}

We assume that the kinematic 3-surface has  constant scalar
curvature $\mathbf{K}\not= 0$. From (\ref{eq311}), we have
\begin{eqnarray*}
P\Big(\phi^{n_1}\cos(m_1\,\phi),\phi^{n_1}\sin(m_1\,\phi)\Big)
-\mathbf{K}\,Q\Big(\phi^{n_2}\cos(m_2\,\phi),\phi^{n_2}\sin(m_2\,\phi)\Big)\\
=\sum_{j=0}^{6}\sum_{i=0}^{6}\Big(A_{i,j}\phi^{i}\cos(j\,\phi)
+B_{i,j}\phi^i\sin(j\,\phi)\Big)=0.
\end{eqnarray*}
In this case, a straightforward computation shows that the
coefficients of $\phi^6\cos(6\phi)$ and $\phi^6\sin(6\phi)$ are
\begin{eqnarray*}
A_{6,6}&=&\frac{1}{16}\,\lambda^6\,\mathbf{K}\,\Big(
\omega_7^6-15\omega_7^4\,\omega_2^2+15\omega_7^2\,\omega_2^2-\omega_2^6\Big),\nonumber\\
B_{6,6}&=&-\frac{1}{8}\,\lambda^6\,\mathbf{K}\,\omega_2\,\omega_7\,\Big(
3\omega_7^4-10\omega_7^2\,\omega_2^2+3\omega_2^4\Big). \nonumber
\end{eqnarray*}
By solving the two equations  $A_{6,6}=0$ and $B_{6,6}=0$, we get
$\omega_2=\omega_7=0$. New computations yields
$$A_{6,0}=-2\,\mathbf{K}\,\Big(\alpha_{2}(1+\lambda^2)-\lambda^4\,s^{\prime2}\Big)$$
From $A_{6,0}=0$ and since both terms $\mathbf{K}$ and
$\Big(\alpha_{2}(1+\lambda^2)-\lambda^4\,s^{\prime2}\Big)$ do not
vanish (see (\ref{a02}) again), we arrive a contradiction.

As a consequence of these computations, we conclude the next result,
which was announced in the Introduction of this work.

\begin{thm}\label{th-32} There are not kinematic 3-surface in $\mathbf{E}^7$ generated by the equifom
motion of a helix whose scalar curvature $\mathbf{K}$ is a non-zero
constant.
\end{thm}


\end{document}